\documentclass[12pt]{article}
\usepackage{amsmath,amsthm,amssymb}
\usepackage{amssymb,latexsym}

\newtheorem{theorem}{Theorem}

\begin{document}
\baselineskip=17pt

\title{\bf The ternary Goldbach problem with prime numbers of a mixed type}

\author{\bf S. I. Dimitrov}

\date{\bf2017}
\maketitle

\begin{abstract}
In the present  paper we prove that every sufficiently large odd integer $N$ can be represented in the form
\begin{equation*}
N=p_1+p_2+p_3\,,
\end{equation*}
where $p_1,p_2,p_3$ are primes, such that $p_1=x^2 + y^2 +1$,\, $p_2=[n^c]$.
\medskip

\textbf{Keywords}: Goldbach problem, Prime numbers, Circle method.
\end{abstract}

\section{Notations}
\indent

Let $N$ be a sufficiently large odd integer.
The letter $p$, with or without subscript, will always denote prime numbers.
Let $A>100$ be a constant.
By $\varepsilon$ we denote an arbitrary small positive number,
not the same in all appearances.
The relation $f(x)\ll g(x)$ means that $f(x)=\mathcal{O}(g(x))$.
As usual $[t]$ and $\{t\}$ denote the integer part, respectively, the
fractional part of $t$. Instead of $m\equiv n\,\pmod {k}$ we write for simplicity $m\equiv n\,(k)$.
As usual $e(t)$=exp($2\pi it$). We denote by $(d,q)$, $[d,q]$ the greatest
common divisor and the least common multiple of $d$ and $q$ respectively.
As usual $\varphi(d)$ is Euler's function; $\mu (d)$ is M\"{o}bius' function;
$r(d)$ is the number of solutions of the
equation $d=m_1^2+m_2^2$ in integers $m_j$; $\chi(d)$ is the non-principal character modulo 4
and $L(s,\chi)$ is the corresponding Dirichlet's $L$-function.
By $c_0$ we denote some positive number, not necessarily the same in different occurrences.
Let $c$ be a real constant such that $1<c<73/64$.

Denote
\begin{align}
\label{gamma}
&\gamma=1/c\,;\\
\label{D}
&D=\frac{N^{1/2}}{(\log N)^A}\,;\\
\label{psi}
&\psi(t)=\{t\}-1/2\,;\\
\label{theta0}
&\theta_0=\frac{1}{2}-\frac{1}{4}e\log2=0.0289...\,;\\
\label{Sigmadl}
&\mathfrak{S}_{d,l}(N)=\prod\limits_{p\nmid d\atop{p|N}}\left(1-\frac{1}{(p-1)^2}\right)
\prod\limits_{p|d\atop{p\nmid N-l}}\left(1-\frac{1}{(p-1)^2}\right)\nonumber\\
&\quad\quad\quad\quad\;\,\times\prod\limits_{p\nmid dN}\left(1+\frac{1}{(p-1)^3}\right)
\prod\limits_{p|d\atop{p|N-l}}\left(1+\frac{1}{p-1}\right)\,;\\
\label{SigmaN}
&\mathfrak{S}(N)=\prod\limits_{p|N}\left(1-\frac{1}{(p-1)^2}\right)
\prod\limits_{p\nmid N}\left(1+\frac{1}{(p-1)^3}\right)\,;\\
\label{SigmaGamma}
&\mathfrak{S}_\Gamma(N)=\pi\mathfrak{S}(N)
\prod\limits_{p\nmid N(N-1)}\left(1+\chi(p)\frac{p-3}{p(p^2-3p+3)}\right)
\prod\limits_{p\mid N} \left(1+\chi(p)\frac{1}{p(p-1)}\right)\nonumber\\
&\quad\quad\quad\quad\times\prod\limits_{p\mid N-1}\left(1+\chi(p)\frac{2p-3}{p(p^2-3p+3)}\right)\,;\\
\label{Delta}
&\Delta(t,h)=\max\limits_{y\leq t}\max\limits_{(l,h)=1}
\left|\sum\limits_{p\leq y\atop{p\equiv l\,(h)}}\log p-\frac{y}{\varphi(h)}\right|\,.
\end{align}

\section{Introduction and statement of the result}
\indent

In 1937 I. M. Vinogradov \cite{Vino} solved the ternary Goldbach problem.
He proved that for a sufficiently large odd integer $N$

\begin{equation*}
\sum\limits_{p_1+p_2+p_3=N}\log p_1\log p_2\log p_3
=\frac{1}{2}\mathfrak{S}(N)N^2
+\mathcal{O}\left(\frac{N^2}{\log^AN}\right)\,,
\end{equation*}
where $\mathfrak{S}(N)$ is defined by \eqref{SigmaN} and $A>0$ is an arbitrarily large constant.

In 1953 Piatetski-Shapiro proved that for any fixed $c\in(1,12/11)$ the sequence
\begin{equation*}
([n^c])_{n\in \mathbb{N}}
\end{equation*}
contains infinitely many prime numbers.
Such prime numbers are named in honor of Piatetski-Shapiro.
The interval for $c$ was subsequently improved many times and the best result
up to now belongs to Rivat and Wu \cite{Rivat-Wu} for $c\in(1,243/205)$.

In 1992, A. Balog and J. P. Friedlander \cite{Balog}  considered the ternary Goldbach problem
with variables restricted to Piatetski-Shapiro primes. They proved that, for
any fixed $1<c<21/20$, every sufficiently large odd integer $N$ can be represented in the form
\begin{equation*}
N=p_1+p_2+p_3\,,
\end{equation*}
where $p_1,p_2,p_3$ are primes, such that $p_k=[n^c_k]$, k=1,2,3.
Rivat \cite{Rivat-Wu} extended the range to $1<c<199/188$; Kumchev \cite{Kumchev}
extended the range to $1<c<53/50$. Jia \cite{Jia} used a sieve method
to enlarge the range to $1<c<16/15$.
Furthermore Kumchev \cite{Kumchev} proved that for any fixed $1<c<73/64$
every sufficiently large odd integer may be written as the sum of two primes
and  prime number of type $p=[n^c]$.

On the other hand in 1960 Linnik \cite{Linnik} showed that there exist infinitely many prime numbers of the form
$p=x^2 + y^2 +1$, where $x$ and $y$ are integers.
In 2010 Tolev \cite{T2} proved that every sufficiently large odd integer $N$ can be represented in the form
\begin{equation*}
N=p_1+p_2+p_3\,,
\end{equation*}
where $p_1,p_2,p_3$ are primes, such that $p_k=x_k^2 + y_k^2 +1$, k=1,2.
In 2017  Ter\"{a}v\"{a}inen \cite{Teravainen} improved Tolev's result
for primes $p_1,p_2,p_3$, such that $p_k=x_k^2 + y_k^2 +1$, k=1,2,3.

Recently the author \cite{Dimitrov} proved that  there exist infinitely many
arithmetic progressions of three different primes $p_1,p_2,p_3=2p_2-p_1$
such that $p_1=x^2 + y^2 +1$,\; $p_3=[n^c]$.

Define
\begin{equation}\label{Gamma}
\Gamma(N)=\sum\limits_{p_1+p_2+p_3=N\atop{p_2=[n^c]}}r(p_1-1)
p^{1-\gamma}_2\log p_1\log p_2\log p_3\,.
\end{equation}
Motivated by these results  we shall prove the following theorem.
\begin{theorem}\label{Golbach-Dimitrov}Assume that $1<c<73/64$. Then the asymptotic formula
\begin{equation*}
\Gamma(N)=\frac{\gamma}{2}\mathfrak{S}_\Gamma(N)N^2
+\mathcal{O}\big(N^2(\log N)^{-\theta_0}(\log\log N)^6\big)\,,
\end{equation*}
holds. Here $\gamma$, $\theta_0$ and $\mathfrak{S}_\Gamma(N)$ are defined by
\eqref{gamma}, \eqref{theta0} and \eqref{SigmaGamma}.
\end{theorem}
Bearing in mind that $\mathfrak{S}_\Gamma(N)\gg1$ for $N$ odd, from Theorem 1 it follows that for any fixed $1<c<73/64$
every sufficiently large odd integer $N$ can be written  in the form
\begin{equation*}
N=p_1+p_2+p_3\,,
\end{equation*}
where $p_1,p_2,p_3$ are primes, such that $p_1=x^2 + y^2 +1$,\, $p_2=[n^c]$.

The asymptotic formula obtained for $\Gamma(N)$ is the product of the individual asymptotic formulas
\begin{equation*}
\sum\limits_{p_1+p_2+p_3=N}r(p_1-1)\log p_1\log p_2\log p_3\sim\frac{1}{2}\mathfrak{S}_\Gamma(N)N^2
\end{equation*}
and
\begin{equation*}
\frac{1}{N}\sum\limits_{p\leq N\atop{p=[n^c]}}p^{1-\gamma}\log p\sim\gamma\,.
\end{equation*}
The proof of Theorem \ref{Golbach-Dimitrov} follows the same ideas as the proof in \cite{Dimitrov}.

\section{Outline of the proof}
\indent

Using \eqref{Gamma} and well-known identity $r(n)=4\sum_{d|n}\chi(d)$ we find
\begin{equation} \label{Gamm}
\Gamma(N)=4\big(\Gamma_1(N)+\Gamma_2(N)+\Gamma_3(N)\big),
\end{equation}
where
\begin{align}
\label{Gamma1}
&\Gamma_1(N)=\sum\limits_{p_1+p_2+p_3=N\atop{p_2=[n^c]}}\left(\sum\limits_{d|p_1-1
\atop{d\leq D}}\chi(d)\right)p^{1-\gamma}_2\log p_1\log p_2\log p_3\,,\\
\label{Gamma2}
&\Gamma_2(N)=\sum\limits_{p_1+p_2+p_3=N\atop{p_2=[n^c]}}\left(\sum\limits_{d|p_1-1
\atop{D<d<N/D}}\chi(d)\right)p^{1-\gamma}_2\log p_1\log p_2\log p_3\,,\\
\label{Gamma3}
&\Gamma_3(N)=\sum\limits_{p_1+p_2+p_3=N\atop{p_2=[n^c]}}\left(\sum\limits_{d|p_1-1
\atop{d\geq N/D}}\chi(d)\right)p^{1-\gamma}_2\log p_1\log p_2\log p_3\,.
\end{align}
In order to estimate $\Gamma_1(N)$ and $\Gamma_3(N)$ we have to consider
the sum
\begin{equation}\label{Ild}
I_{d,l;J}(N)=\sum\limits_{p_1+p_2+p_3=N\atop{p_1\equiv l\,(d)\atop{p_1\in J
\atop{p_2=[n^c]}}}}p^{1-\gamma}_2\log p_1\log p_2\log p_3\,,
\end{equation}
where $d$ and $l$ are coprime natural numbers, and $J\subset[1,N]$.
If $J=[1,N]$ then we write for simplicity $I_{d,l}(N)$.
We apply the circle method.
Clearly
\begin{equation}\label{IldInt}
I_{d,l;J}(N)=\int\limits_{0}^{1}S_{d,l;J}(\alpha)S(\alpha)S_c(\alpha)e(-N\alpha)d\alpha\,,
\end{equation}
where
\begin{align}
\label{Sldj}
&S_{d,l;J}(\alpha)=\sum\limits_{p\in J\atop{p\equiv l\,(d)}}e(\alpha p)\log p\,,\\
\label{S}
&S(\alpha)=S_{1,1;[1,N]}(\alpha)\,,\\
\label{Sc}
&S_c(\alpha)=\sum\limits_{p\leq N\atop{p=[n^c]}}p^{1-\gamma}e(\alpha p)\log p\,.
\end{align}
We define major and minor arcs by
\begin{equation}\label{Arcs}
E_1=\bigcup_{q \le Q} \bigcup_{\substack{a=0 \\ (a, q)=1 }}^{q-1}
 \left[ \frac{a}{q} - \frac{1}{q\tau}, \frac{a}{q} + \frac{1}{q\tau} \right],\,\,\,
E_2= \left[\frac{1}{\tau},1+\frac{1}{\tau} \right]\setminus E_1\,,
\end{equation}
where
\begin{equation}\label{Qtau}
Q=(\log N)^B\,,\;\;\tau=NQ^{-1}\,,\;\;A>4B+3\,,\;\;B>14\,.
\end{equation}
Then we have the decomposition
\begin{equation}\label{Ilddecomp}
I_{d,l;J}(N)=I_{d,l;J}^{(1)}(N)+I_{d,l;J}^{(2)}(N)\,,
\end{equation}
where
\begin{equation}\label{I1I2}
I_{d,l;J}^{(i)}(N)=\int\limits_{E_i}S_{d,l;J}(\alpha)S(\alpha)S_c(\alpha)e(-N\alpha)d\alpha\,,\;\;\;i=1,2.
\end{equation}
We shall estimate $I_{d,l;J}^{(1)}(N)$, $\Gamma_3(N)$, $\Gamma_2(N)$ and $\Gamma_1(N)$
respectively in the sections 4, 5, 6 and 7.
In section 8 we shall complete the proof of the Theorem.
\section{Asymptotic formula for $\mathbf{I_{d,l;J}^{(1)}(N)}$}
\indent

We have
\begin{equation}\label{IldH}
I_{d,l;J}^{(1)}(N)=\sum\limits_{q\le Q}\sum\limits_{\substack{a=0 \\ (a, q)=1 }}^{q-1}H(a,q)\,,
\end{equation}
where
\begin{equation}\label{Haq}
H(a,q)=\int\limits_{-1/q\tau}^{1/q\tau}S_{d,l;J}\left(\frac{a}{q}+\alpha\right)
S\left(\frac{a}{q}+\alpha\right)S_c\left(\frac{a}{q}+\alpha\right)
e\left(-N\left(\frac{a}{q}+\alpha\right)\right)d\alpha\,.
\end{equation}
On the other hand
\begin{equation}\label{SldJT}
S_{d,l;J}\left(\frac{a}{q}+\alpha\right)=\sum\limits_{1\leq m\leq q\atop{(m,q)=1
\atop{m\equiv l\,((d,q))}}}e\left(\frac{am}{q}\right)T(\alpha)+\mathcal{O}\big(q\log N\big)\,,
\end{equation}
where
\begin{equation*}
T(\alpha)=\sum\limits_{p\in J\atop{p\equiv l\,(d)\atop{p\equiv m\,(q)}}}e(\alpha p)\log p\,.
\end{equation*}
According to Chinese remainder theorem there exists integer $f=f(l,m,d,q)$ such that $(f,[d,q])=1$ and
\begin{equation*}
T(\alpha)=\sum\limits_{p\in J\atop{p\equiv f\,([d,q])}}e(\alpha p)\log p\,.
\end{equation*}
Applying Abel's transformation we obtain
\begin{align}\label{TAbel}
T(\alpha)&=-\int\limits_{J_1}^{J_2}\left(\sum\limits_{J_1<p<t\atop{p\equiv f\,([d,q])}}\log p\right)
\frac{d}{dt}(e(\alpha t))dt+\left(\sum\limits_{p\in J\atop{p\equiv f\,([d,q])}}\log p\right)e(\alpha J_2)\nonumber\\
&=-\int\limits_{J_1}^{J_2}\left(\frac{t-J_1}{\varphi([d,q])}+\mathcal{O}\big(\Delta(J_2,[d,q])\big)\right)
\frac{d}{dt}(e(\alpha t))dt\nonumber\\
&+\left(\frac{J_2-J_1}{\varphi([d,q])}+\mathcal{O}\big(\Delta(J_2,[d,q])\big)\right)e(\alpha J_2)\nonumber\\
&=\frac{1}{\varphi([d,q])}\int\limits_{J_1}^{J_2}e(\alpha t)dt+
\mathcal{O}\big((1+|\alpha|(J_2-J_1))\Delta(J_2,[d,q])\big)\,.
\end{align}
We use the well known formula
\begin{equation}\label{SumInt}
\int\limits_{J_1}^{J_2}e(\alpha t)dt=M_J(\alpha)+\mathcal{O}(1)\,,
\end{equation}
where
\begin{equation*}
M_J(\alpha)=\sum\limits_{m\in J}e(\alpha m)\,.
\end{equation*}
Bearing in mind that $|\alpha|\leq1/q\tau$ and $J\subset(1\,,N]$, from (\ref{Qtau}), (\ref{TAbel}) and (\ref{SumInt})
we get
\begin{equation}\label{Tform}
T(\alpha)=\frac{M_J(\alpha)}{\varphi([d,q])}+\mathcal{O}\bigg(\left(1+\frac{Q}{q}\right)\Delta(N,[d,q])\bigg)\,.
\end{equation}
From (\ref{SldJT}) and (\ref{Tform}) it follows
\begin{equation}\label{SldJcdMJ}
S_{d,l;J}\left(\frac{a}{q}+\alpha\right)=\frac{c_d(a,q,l)}{\varphi([d,q])}M_J(\alpha)+
\mathcal{O}\big(Q(\log N)\Delta(N,[d,q])\big)\,,
\end{equation}
where
\begin{equation*}
c_d(a,q,l)=\sum\limits_{1\leq m\leq q\atop{(m,q)=1\atop{m\equiv l\,((d,q))}}}e\left(\frac{am}{q}\right)\,.
\end{equation*}
We shall find asymptotic formula for $S_c\left(\frac{a}{q}+\alpha\right)$.
From \eqref{Sc} we have
\begin{align}\label{ScOmegaSigma}
S_c(\alpha)&=\sum\limits_{p\leq N}p^{1-\gamma}
\big([-p^\gamma]-[-(p+1)^\gamma]\big)e(\alpha p)\log p\nonumber\\
&=\Omega(\alpha)+\Sigma(\alpha)\,,
\end{align}
where
\begin{align}\label{OmegaSigma1}
&\Omega(\alpha)=\sum\limits_{p\leq N}p^{1-\gamma}
\big((p+1)^\gamma-p^\gamma\big)e(\alpha p)\log p\,,\\
\label{OmegaSigma2}
&\Sigma(\alpha)=\sum\limits_{p\leq N}p^{1-\gamma}
\big(\psi(-(p+1)^\gamma)-\psi(-p^\gamma)\big)e(\alpha p)\log p\,.
\end{align}
According to Kumchev (\cite{Kumchev}, Theorem 2) for $64/73<\gamma<1$ uniformly in $\alpha$ we have that
\begin{equation}\label{Sigmaest1}
\Sigma\left(\frac{a}{q}+\alpha\right)\ll N^{1-\varepsilon}\,.
\end{equation}
On the other hand
\begin{equation}\label{p+1}
(p+1)^\gamma-p^\gamma=\gamma p^{\gamma-1}+\mathcal{O}\left(p^{\gamma-2}\right)\,.
\end{equation}
The formulas \eqref{OmegaSigma1} and \eqref{p+1} give us
\begin{equation}\label{OmegaS}
\Omega(\alpha)=\gamma S(\alpha)+\mathcal{O}(N^\varepsilon)\,,
\end{equation}
where $S(\alpha)$ is defined by \eqref{S}.\\
According to (\cite{Karat}, Lemma 3, \S 10) we have
\begin{equation}\label{S-alphaM}
S\left(\frac{a}{q}+\alpha\right)=\frac{\mu(q)}{\varphi(q)}M(\alpha)
+\mathcal{O}\left(Ne^{-c_0\sqrt{\log N}}\right)\,,
\end{equation}
where
\begin{equation*}
M(\alpha)=\sum\limits_{m\leq N}e(\alpha m)\,.
\end{equation*}
Bearing in mind \eqref{ScOmegaSigma}, \eqref{Sigmaest1}, \eqref{OmegaS} and \eqref{S-alphaM} we obtain
\begin{equation}\label{Scest}
S_c\left(\frac{a}{q}+\alpha\right)=\gamma\frac{\mu(q)}{\varphi(q)}M(\alpha)
+\mathcal{O}\left(Ne^{-c_0\sqrt{\log N}}\right)\,.
\end{equation}
Furthermore, we need the trivial estimates
\begin{equation}\label{Trivestimations}
\left|S_{d,l;J}\left(\frac{a}{q}+\alpha\right)\right|\ll\frac{N\log N}{d}\,,\quad
\left|S\left(\frac{a}{q}+\alpha\right)\right|\ll N\,,\quad
|M(\alpha)|\ll N\,,\quad|\mu(q)|\ll1\,.
\end{equation}
By (\ref{SldJcdMJ}), \eqref{S-alphaM} -- (\ref{Trivestimations})
and the well-known inequality $\varphi(n)\gg n(\log\log n)^{-1}$ we find
\begin{align}\label{Sprod}
&S_{d,l;J}\left(\frac{a}{q}+\alpha\right)S\left(\frac{a}{q}+\alpha\right)
S_c\left(\frac{a}{q}+\alpha\right)e\left(-N\left(\frac{a}{q}+\alpha\right)\right)\nonumber\\
&=\gamma\frac{c_d(a,q,l)\mu^2(q)}{\varphi([d,q])\varphi^2(q)}M_J(\alpha)M^2(\alpha)
e\left(-N\left(\frac{a}{q}+\alpha\right)\right)
+\mathcal{O}\left(\frac{N^3}{d}e^{-c_0\sqrt{\log N}}\right)\nonumber\\
&+\mathcal{O}\left(\frac{N^2Q\log^2N}{q^2}\Delta(N,[d,q])\right)\,.
\end{align}
Having in mind (\ref{Qtau}), (\ref{Haq}) and (\ref{Sprod}) we get
\begin{align}\label{Haqest}
H(a,q)&=\gamma\frac{c_d(a,q,l)\mu^2(q)}{\varphi([d,q])\varphi^2(q)}e\left(-N\frac{a}{q}\right)
\int\limits_{-1/q\tau}^{1/q\tau}M_J(\alpha)M^2(\alpha)e(-N\alpha)d\alpha\nonumber\\
&+\mathcal{O}\left(\frac{N^2}{qd}e^{-c_0\sqrt{\log N}}\right)
+\mathcal{O}\left(\frac{NQ^2\log^2N}{q^3}\Delta(N,[d,q])\right)\,.
\end{align}
Taking into account \eqref{IldH}, \eqref{Haqest}
and following the method in \cite{T1} we obtain
\begin{align}\label{I1}
I_{d,l;J}^{(1)}(N)&=
\gamma\frac{\mathfrak{S}_{d,l}(N)}{\varphi(d)}\sum\limits_{m_1+m_2+m_3=N\atop{m_1\in J}}1
+\mathcal{O}\left(\frac{N^2}{d}(\log N)\sum\limits_{q>Q}\frac{(d,q)\log q}{q^2}\right)\nonumber\\
&+\mathcal{O}\left(\tau^2(\log N)\sum\limits_{q\le Q}\frac{q}{[d,q]}\right)
+\mathcal{O}\left(NQ^2(\log N)^2\sum\limits_{q\leq Q}
\frac{\Delta(N,[d,q])}{q^2}\right)\nonumber\\
&+\mathcal{O}\left(\frac{N^2}{d}e^{-c_0\sqrt{\log N}}\right)\,,
\end{align}
where $\mathfrak{S}_{d,l}(N)$ is defined by \eqref{Sigmadl}.

\newpage

\section{Upper bound for $\mathbf{\Gamma_3(N)}$}
\indent

Consider the sum  $\Gamma_3(N)$.\\
Since
\begin{equation*}
\sum\limits_{d|p_1-1\atop{d\geq N/D}}\chi(d)=\sum\limits_{m|p_1-1\atop{m\leq (p_1-1)D/N}}
\chi\bigg(\frac{p_1-1}{m}\bigg)
=\sum\limits_{j=\pm1}\chi(j)\sum\limits_{m|p_1-1\atop{m\leq (p_1-1)D/N
\atop{\frac{p_1-1}{m}\equiv j\,(4)}}}1
\end{equation*}
then from \eqref{Gamma3} and \eqref{Ild} it follows
\begin{equation*}
\Gamma_3(N)=\sum\limits_{m<D\atop{2|m}}\sum\limits_{j=\pm1}\chi(j)I_{4m,1+jm;J_m}(N)\,,
\end{equation*}
where $J_m=[1+mN/D,N]$.
Therefore from \eqref{Ilddecomp} we get
\begin{equation}\label{Gamm3}
\Gamma_3(N)=\Gamma_3^{(1)}(N)+\Gamma_3^{(2)}(N)\,,
\end{equation}
where
\begin{equation}\label{Gamm3nu}
\Gamma_3^{(\nu)}(N)=\sum\limits_{m<D\atop{2|m}}\sum\limits_{j=\pm1}\chi(j)
I_{4m,1+jm;J_m}^{(\nu)}(N)\,,\;\;\;\nu=1,2.
\end{equation}

Let us consider first $\Gamma_3^{(2)}(N)$. Bearing in mind \eqref{I1I2} for $i=2$ and
\eqref{Gamm3nu} for $\nu=2$ we have
\begin{equation*}
\Gamma_3^{(2)}(N)=\int\limits_{E_2}K(\alpha)S(\alpha)S_c(\alpha)e(-N\alpha)d\alpha\,,
\end{equation*}
where
\begin{equation}\label{Kalpha}
K(\alpha)=\sum\limits_{m<D\atop{2|m}}\sum\limits_{j=\pm1}\chi(j)S_{4m,1+jm;J_m}(\alpha)\,.
\end{equation}
Using Cauchy's inequality we obtain
\begin{align}\label{Cauchy}
\Gamma_3^{(2)}(N)&\ll\sup\limits_{\alpha\in E_2\setminus\{1\}}|S_c(\alpha)
|\int\limits_{E_2}|K(\alpha)S(\alpha)|d\alpha+\mathcal{O}(N^\varepsilon)\nonumber\\
&\ll\sup\limits_{\alpha\in E_2\setminus\{1\}}|S_c(\alpha)|\left(\int\limits_{0}^1|K(\alpha)|^2d\alpha\right)^{1/2}
\left(\int\limits_{0}^1|S(\alpha)|^2d\alpha\right)^{1/2}+\mathcal{O}(N^\varepsilon)\,.
\end{align}
From \eqref{ScOmegaSigma} and \eqref{OmegaS} we have
\begin{align}\label{ScOmegaSigma2}
S_c(\alpha)=\gamma S(\alpha)+\Sigma(\alpha)+\mathcal{O}(N^\varepsilon)\,,
\end{align}
where $S(\alpha)$ and $\Sigma(\alpha)$ are defined by \eqref{S} and \eqref{OmegaSigma2}.\\
Using \eqref{Arcs} and \eqref{Qtau}
we can prove in the same way as in (\cite{Karat}, Ch.10, Th.3) that
\begin{equation} \label{supS}
\sup\limits_{\alpha\in E_2\setminus\{1\}}|S(\alpha)|\ll\frac{N}{(\log N)^{B/2-4}}\,.
\end{equation}
According to Kumchev (\cite{Kumchev}, Theorem 2) we have that
\begin{equation}\label{Sigmaest2}
\sup\limits_{\alpha\in E_2\setminus\{1\}}|\Sigma(\alpha)|\ll N^{1-\varepsilon}\,.
\end{equation}
Bearing in mind \eqref{ScOmegaSigma2} -- \eqref{Sigmaest2} we get
\begin{equation}\label{supSc}
\sup\limits_{\alpha\in E_2\setminus\{1\}}|S_c(\alpha)|\ll\frac{N}{(\log N)^{B/2-4}}\,.
\end{equation}
From \eqref{S} after straightforward computations we find
\begin{equation}\label{IntS}
\int\limits_{0}^1|S(\alpha)|^2d\alpha\ll N\log N\,.
\end{equation}
On the other hand from \eqref{Sldj} and \eqref{Kalpha} we obtain
\begin{align}\label{IntK}
\int\limits_{0}^1|K(\alpha)|^2\,d\alpha&=
\sum\limits_{m_1,m_2<D\atop{2\mid m_1,2\mid m_2}}\sum\limits_{j_1=\pm1\atop{j_2=\pm1}}\chi(j_1)\chi(j_2)\nonumber\\
&\times\int\limits_0^1S_{4m_1, 1+j_1m_1;J_{m_1}}(\alpha)S_{4m_2, 1+j_2m_2;J_{m_2}}(-\alpha)d\alpha\nonumber\\
&=\sum\limits_{m_1,m_2<D\atop{2\mid m_1,2\mid m_2}}\sum\limits_{j_1=\pm1\atop{j_2=\pm1}}\chi(j_1)\chi(j_2)\nonumber\\
&\times\sum\limits_{p_i\in J_{m_i},i=1,2\atop{p_i\equiv 1+j_im_i\,(4m_i),i=1,2}}
\log p_1\log p_2\int\limits_0^1e(\alpha(p_1-p_2))d\alpha\nonumber\\
&=\sum\limits_{m<D\atop{2\mid m}}\sum\limits_{j=\pm1}\chi(j)\sum\limits_{p\in J_m\atop{p\equiv 1+jm\,(4m)}}
(\log p)^2\nonumber\\
&\ll(\log N)^2\sum\limits_{m<D\atop{2\mid m}}\sum\limits_{p\in J_m\atop{p\equiv 1+jm\,(4m)}}1\nonumber\\
&\ll N(\log N)^2\sum\limits_{m<D}\frac{1}{m}\nonumber\\
&\ll N\log^3 N\,.
\end{align}
Thus from \eqref{Cauchy}, \eqref{supSc} -- \eqref{IntK} it follows
\begin{equation} \label{Gamma32est}
\Gamma_3^{(2)}(N)\ll\frac{N^2}{(\log N)^{B/2-6}}\,.
\end{equation}

Now let us consider $\Gamma_3^{(1)}(N)$. From \eqref{I1} and \eqref{Gamm3nu}
for $\nu=1$ we get
\begin{align}
\Gamma_3^{(1)}(N)=\Gamma^*&+\mathcal{O}\big(N^2(\log N)\Sigma_1\big)
+\mathcal{O}\big(\tau^2(\log N)\Sigma_2\big)\nonumber\\
\label{Gamma31}
&+\mathcal{O}\big(NQ^2(\log N)^2\Sigma_3\big)+\mathcal{O}\big(N^2e^{-c_0\sqrt{\log N}}\Sigma_4\big)\,,
\end{align}
where
\begin{align*}
&\Gamma^*=\gamma\Bigg(\sum\limits_{m_1+m_2+m_3=N\atop{m_1\in J_m}}1\Bigg)
\sum\limits_{m<D\atop{2|m}}\frac{1}{\varphi(4m)}\sum\limits_{j=\pm1}\chi(j)\mathfrak{S}_{4m,1+jm}(N)\,,\\
&\Sigma_1=\sum\limits_{m<D}\sum\limits_{q>Q}\frac{(4m,q)\log q}{mq^2}\,,\\
&\Sigma_2=\sum\limits_{m<D}\sum\limits_{q\leq Q}\frac{q}{[4m,q]}\,,\\
&\Sigma_3=\sum\limits_{m<D}\sum\limits_{q\leq Q}\frac{\Delta(N,[4m,q])}{q^2}\,,\\
&\Sigma_4=\sum\limits_{m<D}\frac{1}{m}\,.
\end{align*}
From the definition \eqref{Sigmadl} it follows that $\mathfrak{S}_{4m,1+jm}(N)$ does not depend on $j$.
Then we have $\sum\limits_{j=\pm1}\chi(j)\mathfrak{S}_{4m,1+jm}(N)=0$ and that leads to
\begin{equation} \label{Gamma*}
\Gamma^*=0\,.
\end{equation}
Arguing as in \cite{T1} and using Bombieri -- Vinogradov's theorem
we find the following estimates
\begin{equation}\label{errorterms1,2}
\Sigma_1\ll\frac{\log^3N}{Q}\,,\;\;\;\Sigma_2\ll Q\log^2N\,,
\end{equation}
\begin{equation}\label{errorterms3,4}
\;\;\;\;\;\;\Sigma_3\ll\frac{N}{(\log N)^{A-B-5}}\,,\;\;\;\Sigma_4\ll\log N\,.
\end{equation}
Bearing in mind \eqref{Qtau}, \eqref{Gamma31} -- \eqref{errorterms3,4} we obtain
\begin{equation} \label{Gamma31est}
\Gamma_3^{(1)}(N)\ll\frac{N^2}{(\log N)^{B-4}}\,.
\end{equation}

Now from \eqref{Gamm3}, \eqref{Gamma32est} and \eqref{Gamma31est} we find
\begin{equation}\label{Gamm3est}
\Gamma_3(N)\ll\frac{N^2}{(\log N)^{B/2-6}}\,.
\end{equation}

\newpage

\section{Upper bound for $\mathbf{\Gamma_2(N)}$}
\indent

Consider the sum $\Gamma_2(N)$ defined by \eqref{Gamma2}.\\
We denote by $\mathcal{F}$ the set of all primes
$p\leq N$ such that $p-1$ has a divisor belongs to the interval $(D,N/D)$.
Using the inequality $uv\leq u^2+v^2$ and taking into account the symmetry with respect to $d$ and $t$ we get
\begin{align}\label{Gamma2est1}
\Gamma_2(N)^2&\ll(\log N)^6N^{2-2\gamma}\sum\limits_{p_1+p_2+p_3=N
\atop{p_4+p_5+p_6=N\atop{p_2=[n_1^c],\, p_5=[n_2^c]}}}\bigg|\sum\limits_{d|p_1-1\atop{D<d<N/D}}\chi(d)\bigg|
\bigg|\sum\limits_{t|p_4-1\atop{D<t<N/D}}\chi(t)\bigg|\nonumber\\
&\ll(\log N)^6N^{2-2\gamma}\sum\limits_{p_1+p_2+p_3=N
\atop{p_4+p_5+p_6=N\atop{p_2=[n_1^c],\, p_5=[n_2^c]\atop{p_4\in\mathcal{F}}}}}\bigg|\sum\limits_{d|p_1-1
\atop{D<d<N/D}}\chi(d)\bigg|^2\,.
\end{align}
Further we use that if $n$ is a natural such that $n\leq N$, then the
number of solutions of the equation $p_1+p_2=n$ in primes $p_1,p_2\leq N$
such that $p_1=[m^{1/\gamma}]$ is $\mathcal{O}\big(N^\gamma(\log N)^{-2}\log\log N\big)$, i.e.
\begin{equation}\label{p1p2nN}
\#\{p_1\,:p_1+p_2=n,\;\;p_1=[m^{1/\gamma}],\;\;n\leq N\}
\ll\frac{N^\gamma\log\log N}{\log^2N}\,.
\end{equation}
This follows for example from (\cite{Halberstam}, Ch.2, Th.2.4).\\
Thus the summands in the sum \eqref{Gamma2est1} for which $p_1=p_4$ can be estimated with
$\mathcal{O}(N^{3+\varepsilon})$.\\
Therefore
\begin{equation}\label{Gamma2est2}
\Gamma_2(N)^2\ll(\log N)^6N^{2-2\gamma}\Sigma_1+N^{3+\varepsilon}\,,
\end{equation}
where
\begin{equation*}
\Sigma_1=\sum\limits_{p_1\leq N}\bigg|\sum\limits_{d|p_1-1
\atop{D<d<N/D}}\chi(d)\bigg|^2\sum\limits_{p_4\leq N\atop{p_4\in\mathcal{F}
\atop{p_4\neq p_1}}}\sum\limits_{p_2+p_3=N-p_1\atop{p_5+p_6=N-p_4\atop{p_2=[n_1^c],\, p_5=[n_2^c]}}}1\,.
\end{equation*}
We use again \eqref{p1p2nN} and find
\begin{equation}\label{Sigmaest}
\Sigma_1\ll\frac{N^{2\gamma}}{\log^4N}(\log\log N)^2\Sigma_2\Sigma_3\,,
\end{equation}
where
\begin{equation*}
\Sigma_2=\sum\limits_{p\leq N}\Bigg|\sum\limits_{d|p-1
\atop{D<d<N/D}}\chi(d)\Bigg|^2\,,\;\;\;\;\;
\Sigma_3=\sum\limits_{p\leq N\atop{p\in\mathcal{F}}}1\,.
\end{equation*}
Arguing as in (\cite{Hooley}, Ch.5) we obtain
\begin{equation}\label{Sigma12}
\Sigma_2\ll\frac{N(\log\log N)^7}{\log N}\,,\;\;\;\;\;
\Sigma_3\ll\frac{N(\log\log N)^3}{(\log N)^{1+2\theta_0}}\,.
\end{equation}
where $\theta_0$ is denoted by  \eqref{theta0}.

From \eqref{Gamma2est2} -- \eqref{Sigma12} it follows
\begin{equation}\label{Gamm2est}
\Gamma_2(N)\ll N^2(\log N)^{-\theta_0}(\log\log N)^6\,.
\end{equation}
\section{Asymptotic formula for $\mathbf{\Gamma_1(N)}$}
\indent

In this section our argument is a modification of Tolev's \cite{T2} argument.\\
Consider the sum $\Gamma_1(N)$.
From \eqref{Gamma1}, \eqref{Ild} and \eqref{Ilddecomp} we get
\begin{equation}\label{Gamm1}
\Gamma_1(N)=\Gamma_1^{(1)}(N)+\Gamma_1^{(2)}(N)\,,
\end{equation}
where
\begin{equation*}
\Gamma_1^{(1)}(N)=\sum\limits_{d\leq D}\chi(d)I_{d,1}^{(1)}(N)\,,
\end{equation*}
\begin{equation*}
\Gamma_1^{(2)}(N)=\sum\limits_{d\leq D}\chi(d)I_{d,1}^{(2)}(N)\,.
\end{equation*}
We estimate the sum $\Gamma_1^{(2)}(N)$ by the same way as the
sum $\Gamma_3^{(2)}(N)$ and obtain
\begin{equation}\label{Gamma12est}
\Gamma_1^{(2)}(N)\ll\frac{N^2}{(\log N)^{B/2-6}}\,.
\end{equation}
Now we consider $\Gamma_1^{(1)}(N)$. We use the formula \eqref{I1} for $J=[1,N]$.
The error term is estimated by the same way as for $\Gamma_3^{(1)}(N)$.
We have
\begin{equation}\label{Gamma11est1}
\Gamma_1^{(1)}(N)=\frac{\gamma}{2}\mathfrak{S}(N)N^2
\sum\limits_{d\leq D}\frac{\chi(d)\mathfrak{S}^\ast_{d,1}(N)}{\varphi(d)}
+\mathcal{O}\bigg(\frac{N^2}{(\log X)^{B-4}}\bigg)\,,
\end{equation}
where $\mathfrak{S}(N)$ is defined by \eqref{SigmaN} and
\begin{align}\label{Sigmastar}
&\mathfrak{S}^\ast_{d,1}(N)=\prod\limits_{p|d\atop{p|N}}\left(1-\frac{1}{(p-1)^2}\right)^{-1}
\prod\limits_{p|d\atop{p\nmid N-1}}\left(1-\frac{1}{(p-1)^2}\right)\nonumber\\
&\quad\quad\quad\quad\;\,\times\prod\limits_{p|d\atop{p\nmid N}}\left(1+\frac{1}{(p-1)^3}\right)^{-1}
\prod\limits_{p|d\atop{p|N-1}}\left(1+\frac{1}{p-1}\right)\,;
\end{align}
Denote
\begin{equation}\label{Sigmaf}
\Sigma=\sum\limits_{d\leq D}f(d)\,,\;\;\;\;\;f(d)=\frac{\chi(d)\mathfrak{S}^\ast_{d,1}(N)}{\varphi(d)}\,.
\end{equation}
We have
\begin{equation}\label{fdest}
f(d)\ll d^{-1}(\log\log(10d))^2
\end{equation}
with absolute constant in the Vinogradov's symbol. Hence the corresponding Dirichlet series
\begin{equation*}
F(s)=\sum\limits_{d=1}^\infty\frac{f(d)}{d^s}
\end{equation*}
is absolutely convergent in $Re(s)>0$. On the other hand $f(d)$ is a multiplicative with
respect to $d$ and applying Euler's identity we find
\begin{equation}\label{FT}
F(s)=\prod\limits_pT(p,s)\,,\;\;\;\;T(p,s)=1+\sum\limits_{l=1}^\infty f(p^l)p^{-ls}\,.
\end{equation}
From \eqref{Sigmastar}, (\ref{Sigmaf}) and (\ref{FT}) we establish that
\begin{equation*}
T(p,s)=\left(1-\frac{\chi(p)}{p^{s+1}}\right)^{-1}\left(1+\frac{\chi(p)}{p^{s+1}}E_d(p)\right)\,,
\end{equation*}
where
\begin{equation*}
E_d(p)=\begin{cases}(p-3)(p^2-3p+3)^{-1}\;\;\;\;\;\mbox{if}\;\; p\nmid N(N-1)\,,\\
(p-1)^{-1}\;\;\; \quad \quad\quad \quad\quad \quad \;\mbox{if}\;\; p\mid N\,,
\\ (2p-3)(p^2-3p+3)^{-1}\;\;\;\mbox{if}\;\; p\mid N-1\,.
\end{cases}
\end{equation*}
Hence we find
\begin{equation}\label{Fs}
F(s)=L(s+1,\chi)\mathcal{N}(s)\,,
\end{equation}
where $L(s+1,\chi)$ is Dirichlet series corresponding to the character $\chi$ and
\begin{align}\label{Ns}
\mathcal{N}(s)=&\prod\limits_{p\nmid N(N-1)}\left(1+\chi(p)\frac{p-3}{p^{s+1}(p^2-3p+3)}\right)
\prod\limits_{p\mid N} \left(1+\chi(p)\frac{1}{p^{s+1}(p-1)}\right)\nonumber\\
&\times\prod\limits_{p\mid N-1} \left(1+\chi(p)\frac{2p-3}{p^{s+1}(p^2-3p+3)}\right)\,.
\end{align}
From the properties of the L-functions it follows that $F(s)$ has an analytic continuation to $Re(s)>-1$.
It is well known that
\begin{equation}\label{Lsest}
L(s+1,\chi)\ll1+\left|Im(s)\right|^{1/6}\;\;\;\mbox{for}\;\;Re(s)\geq-\frac{1}{2}\,.
\end{equation}
Moreover
\begin{equation}\label{Nsest}
\mathcal{N}(s)\ll1\,.
\end{equation}
Using (\ref{Fs}), (\ref{Lsest}) and (\ref{Nsest}) we get
\begin{equation}\label{Fsest}
F(s)\ll N^{1/6}\;\;\;\mbox{for}\;\;Re(s)\geq-\frac{1}{2}\,,\;\;|Im(s)|\leq N\,.
\end{equation}
We apply Perron's formula given at Tenenbaum (\cite{Tenenbaum}, Chapter II.2) and also (\ref{fdest}) to
obtain
\begin{equation}\label{SigmaPeron}
\Sigma=\frac{1}{2\pi\imath}\int\limits_{\varkappa-\imath N}^{\varkappa+\imath N}F(s)\frac{D^s}{s}ds
+\mathcal{O}\left(\sum\limits_{t=1}^\infty\frac{D^\varkappa\log\log(10t)}{t^{1+\varkappa}
\left(1+N\left|\log\frac{D}{t}\right|\right)}\right)\,,
\end{equation}
where $\varkappa=1/10$. It is easy to see that the error term above is $\mathcal{O}\left(N^{-1/20}\right)$.
Applying the residue theorem we see that the main term in \eqref{SigmaPeron} is equal to
\begin{equation*}
F(0)+\frac{1}{2\pi\imath}\left(\int\limits_{1/10-\imath N}^{-1/2-\imath N}+
\int\limits_{-1/2-\imath N}^{-1/2+\imath N}
+\int\limits_{-1/2+\imath N}^{1/10+\imath N}\right)F(s)\frac{D^s}{s}ds\,.
\end{equation*}
From (\ref{Fsest}) it follows that the contribution from the above integrals is $\mathcal{O}\left(N^{-1/20}\right)$.\\
Hence
\begin{equation}\label{Sigmaest}
\Sigma=F(0)+\mathcal{O}\left(N^{-1/20}\right)\,.
\end{equation}
Using (\ref{Fs}) we get
\begin{equation}\label{F0}
F(0)=\frac{\pi}{4}\mathcal{N}(0)\,.
\end{equation}
Bearing in mind (\ref{Gamma11est1}), (\ref{Sigmaf}), (\ref{Ns}), (\ref{Sigmaest})
and (\ref{F0}) we find a new expression for $\Gamma_1^{(1)}(N)$
\begin{equation}\label{Gamma11est2}
\Gamma_1^{(1)}(N)=\frac{\gamma}{8}\mathfrak{S}_\Gamma(N)N^2
+\mathcal{O}\bigg(\frac{N^2}{(\log N)^{B-4}}\bigg)\,,
\end{equation}
where $\mathfrak{S}_\Gamma$ is defined by (\ref{SigmaGamma}).

From  \eqref{Gamm1}, \eqref{Gamma12est} and \eqref{Gamma11est2} we obtain
\begin{equation} \label{Gamm1est}
\Gamma_1(N)=\frac{\gamma}{8}\mathfrak{S}_\Gamma(N)N^2
+\mathcal{O}\bigg(\frac{N^2}{(\log N)^{B/2-6}}\bigg)\,.
\end{equation}
\section{Proof of the Theorem}
\indent

Therefore using \eqref{Gamm}, \eqref{Gamm3est}, \eqref{Gamm2est} and \eqref{Gamm1est}
we find
\begin{equation*}
\Gamma(N)=\frac{\gamma}{2}\mathfrak{S}_\Gamma(N)N^2
+\mathcal{O}\big(N^2(\log N)^{-\theta_0}(\log\log N)^6\big)\,.
\end{equation*}
This implies that $\Gamma(N)\rightarrow\infty$ as $N\rightarrow\infty$.

The Theorem is proved.

\vskip20pt
\footnotesize
\begin{flushleft}
S. I. Dimitrov\\
Faculty of Applied Mathematics and Informatics\\
Technical University of Sofia \\
8, St.Kliment Ohridski Blvd. \\
1756 Sofia, BULGARIA\\
e-mail: sdimitrov@tu-sofia.bg\\
\end{flushleft}

\begin{thebibliography}{}

\bibitem{Balog}A. Balog, J. P. Friedlander, {\it A hybrid of theorems of Vinogradov and Piatetski-Shapiro},
Pacific J. Math., {\bf156}, (1992), 45 -- 62.

\bibitem{Dimitrov}S. I. Dimitrov, {\it  Prime triples $p_1,p_2,p_3$ in arithmetic progressions
such that $p_1=x^2+y^2+1$, $p_3=[n^c]$}, Notes on Number Theory and Discrete Mathematics, \textbf{23}(4), (2017), 22 -- 33.

\bibitem{Halberstam}H. Halberstam, H.-E. Richert, {\it Sieve Methods},
Academic Press, (1974).

\bibitem{Hooley}C. Hooley, {\it Applications of sieve methods to the theory of numbers},
Cambridge Univ. Press, (1976).

\bibitem{Jia}C.-H. Jia, {\it On the Piatetski-Shapiro-Vinogradov theorem},
Acta Arith., {\bf73}, (1995), 1 -- 28.

\bibitem{Karat}A. Karatsuba, {\it Principles of the Analytic Number Theory},
Nauka, Moscow, (1983) (in Russian).

\bibitem{Kumchev}A. Kumchev, {\it On the Piatetski-Shapiro-Vinogradov Theorem},
Journal de Th\'{e}orie des Nombres de Bordeaux, {\bf9}, (1997), 11 -- 23.

\bibitem{Linnik}Ju. Linnik, {\it An asymptotic formula in an additive problem of Hardy and Littlewood},
Izv. Akad. Nauk SSSR, Ser.Mat., {\bf24}, (1960), 629 -- 706 (in Russian).

\bibitem{Shapiro}I. I. Piatetski-Shapiro, {\it On the distribution of prime numbers in sequences of the form $[f(n)]$},
Mat. Sb., {\bf 33}, (1953), 559 -- 566.

\bibitem{Rivat-Wu}J. Rivat, J. Wu, {\it Prime numbers of the form $[n^c]$},
Glasg. Math. J, {\bf 43}, 2, (2001), 237 -- 254.

\bibitem{Tenenbaum}G. Tenenbaum, {\it Introduction to Analytic and Probabilistic Number Theory},
Cambridge Univ. Press, (1995).

\bibitem{Teravainen}J. Ter\"{a}v\"{a}inen, {\it The Goldbach problem for primes that are sums of two squares plus one},
Mathematika, {\bf64}, 1, (2018), 20 -- 70.

\bibitem{T1}D. Tolev, {\it On the number of representations of an odd integer as a sum of three primes, one of which belongs
to an arithmetic progression},
Proc. Steklov Math. Inst., {\bf 218}, (1997),  415 -- 432.

\bibitem{T2}D. Tolev, {\it The ternary Goldbach problem with arithmetic weights attached to one of the variables},
J. Number Theory, {\bf130}, (2010), 439 -- 457.

\bibitem{Vino}I. M. Vinogradov, {\it Representation of an odd number as the sum of three primes},
Dokl. Akad. Nauk. SSSR, {\bf15}, (1937), 291 -- 294, (in Russian).
\end{thebibliography}
\end{document}